\documentclass[a4paper,10pt]{article}
\usepackage{amsmath,amsthm,amssymb,amscd,epsfig}

\newtheorem{theo}{Theorem}
\newtheorem{lem}[theo]{Lemma}

\newtheorem{coro}[theo]{Corollary}

\def\A{\mathcal A}
\def\R{\mathbb R}
\def\C{\mathbb C}

\def\D{\mathbb D}

\def\H{\mathcal H}
\def\M{\mathcal M}

\def\B{\mathcal B}
\def\diam{\text{diam}}
\def\Lip{\operatorname{Lip}}
\def\lip{\operatorname{lip}}

\title{Removable singularities for H\"older continuous quasiregular mappings in the plane}

\author{Albert Clop
\thanks{The author was supported by projects MTM2004-00519, HF2004-0208, 2001-SGR-00431
\newline Keywords   Quasiconformal, Hausdorff measure, Removability}
}
\date{}

\begin{document}

\maketitle
\begin{abstract}
We give necessary conditions for a set $E$ to be removable for H\"older continuous quasiregular mappings in the plane. We also obtain some removability results for H\"older continuous mappings of finite distortion. 
\end{abstract}

\section{Introduction}

Let $\alpha\in(0,1]$. A function $f:\C\rightarrow\C$ is said to be $\alpha$-H\"older continuous, that is, $f\in\Lip_\alpha(\C)$, if
\begin{equation}\label{lipalfa}
|f(z)-f(w)|\leq C\,|z-w|^\alpha
\end{equation}
whenever $z,w\in\C$ are such that $|z-w|<1$. A set $E\subset\C$ is said to be {\it{removable}} for $\alpha$-H\"older continuous analytic functions if every function $f\in\Lip_\alpha(\C)$, holomorphic on $\C\setminus E$, is actually an entire function. It turns out that there is a characterization of these sets $E$ in terms of Hausdorff measures. For $\alpha\in(0,1)$, Dolzenko \cite{Do} proved that a set $E$ is removable for $\alpha$-H\"older continuous analytic functions if and only if $\H^{1+\alpha}(E)=0$. In particular, for any set $E$ with positive $(1+\alpha)$-dimensional Hausdorff measure, there is a function $f\in\Lip_\alpha(\C)$, holomorphic on $\C\setminus E$, which is not entire. For $\alpha=1$, we deal with the class of Lipschitz continuous analytic functions. Although the same characterization holds, a more involved argument, due to Uy \cite{U}, is needed to show that sets of positive area are not removable.\\
We can ask ourselves which is the corresponding situation in the more general setting of $K$-quasiregular mappings. Given a domain $\Omega\subset\C$ and $K\geq 1$, we say that a mapping $f:\Omega\rightarrow\C$ is $K$-quasiregular in $\Omega$ if $f$ is a $W^{1,2}_{loc}(\Omega)$ solution of the Beltrami equation,
$$\overline\partial f(z)=\mu(z)\,\partial f(z)$$
for almost every $z\in\Omega$, where $\mu$, the Beltrami coefficient, is a measurable function such that $|\mu(z)|\leq\frac{K-1}{K+1}$ at almost every $z\in\Omega$. If $f$ is a homeomorphism, then we say that $f$ is $K$-quasiconformal. When $\mu=0$, we recover the classes of analytic functions and conformal mappings on $\Omega$, respectively.\\
\\
The purpose of this paper is to give sufficient conditions for a compact set $E\subset\C$ to be removable for $\Lip_\alpha$ $K$-quasiregular mappings. In other words, which sets $E$ are such that every $\alpha$-H\"older continuous function $f:\C\rightarrow\C$, $K$-quasiregular on $\C\setminus E$, is actually $K$-quasiregular on the whole plane. In fact, $K$-quasiregular mappings may also be defined in $\R^n$, so that this problem makes sense even in higher dimensions. In this paper, we restrict ourselves to the planar case.\\
\\
Some results have already been given related to this subject. For instance, Koskela and Martio \cite{KM} showed that if $d=\frac{1}{K}\left(1+\frac{\alpha}{K}\right)$, then every compact set $E\subset\C$ with $\H^{d}(E)=0$ is removable for $K$-quasiregular mappings in $\Lip_\alpha$. Moreover, they also gave some sufficient conditions for removability in $\R^n$. In particular, they proved that if $\H^\lambda(E)=0$, $\lambda=\min\{1,n\alpha\}$, then $E$ is removable for $K$-quasiregular mappings in $\Lip_\alpha$. On the other hand, while dealing with 
second order quasilinear elliptic equations, Kilpel\"ainen and Zhong \cite{KZ} showed that $\H^{\alpha(n-1)}(E)=0$ is enough. In this paper, we show that compact sets $E\subset\C$ of zero $d$-dimensional Hausdorff measure, $d=\frac{2}{K+1}(1+\alpha K)$, are removable for $\alpha$-H\"older continuous planar $K$-quasiregular mappings. This index was suggested to us by K. Astala.\\
\\
Unfortunately, our method does not work in $\R^n$, since there is no factorization theorem when $n>2$. However, our result may be rewritten in the planar class of finite distortion mappings. Recall that $f:\Omega\subset\C\rightarrow\C$ is said to be a mapping of finite distortion if $f\in W^{1,1}_{loc}(\Omega)$, $Jf\in L^1_{loc}(\Omega)$ and
$$|Df(x)|^2\leq K(x)\,Jf(x)$$
at almost every $z\in\Omega$. Here $K$ is a measurable function, $K(z)\in (1,\infty)$ at almost every $z\in\Omega$. When $n=2$, this is equivalent to $|\mu(z)|<1$ for almost every $z\in\Omega$. If $K\in L^\infty$, we get the class of bounded distortion (or quasiregular) functions. During the last years, there has been deep progresses in the knowledge of these mappings. Questions like the removability of singularities for bounded mappings of finite distortion have been studied (see, for instance, \cite{AIKM}, \cite{FKZ}, \cite{IKMS} or \cite{IM}). In this planar setting, relatively weak assumptions on $\mu$ are sufficient to ensure that some basic properties of the quasiregular case also hold. For instance, continuity, discreteness and openness, as soon as the existence of normalized solutions and factorization theorems. Althoug the natural regularity for these mappings is the Orlicz-Sobolev space $W^{1,P}_{loc}(\C)$, $P(t)=\frac{t^2}{\log(e+t)}$ (which is larger than the usual $W^{1,2}_{loc}(\C)$), we can study the removability problem for $\Lip_\alpha$ mappings of finite distortion (see Corollaries 3 and 4), which appears as a limiting situation of the quasiregular case.

\section{Sufficient conditions for $\Lip_\alpha$ removability}
We say that a $K$-quasiconformal mapping $\phi:\C\rightarrow\C$ is {\it{principal}} if $\phi(z)-z=O(1/|z|)$ as $|z|\to\infty$. Then, for compactly supported $\mu$ there exist exactly one principal $K$-quasiconformal mapping $\phi$ such that $\overline\partial\phi=\mu\,\partial\phi$ and $\phi(0)=0$, which we will refer to as the {\it{principal solution}} of this Beltrami equation. We denote by $A$ the planar Lebesgue measure, while $\H^d$ stands for Hausdorff $d$-dimensional measure. Given a disk $D$ of radius $r$, we will denote by $\lambda D$ the disk concentric with $D$ whose radius is $\lambda r$.\\
\\We start with an auxiliary result.

\begin{lem}\label{convergence}
Suppose that $\mu_n,\mu$ are measurable functions, compactly supported on $\D$, and assume that $\|\mu_n\|_\infty,\|\mu\|_\infty\leq\frac{K-1}{K+1}$. Let $\phi_n,\phi:\C\rightarrow\C$ be, respectively, the principal solutions to the corresponding Beltrami 
equations. If $\mu_n\to\mu$ at almost every point, then:
\begin{enumerate}
\item $\displaystyle\lim_n J\phi_n=J\phi$, with convergence in $L^p_{loc}$, for every $p\in(1,\frac{K}{K-1})$.
\item $\displaystyle\lim_n \phi_n=\phi$ uniformly on compact sets.
\item $\displaystyle\lim_n \phi_n^{-1}=\phi^{-1}$ uniformly on compact sets.
\end{enumerate}
\end{lem}
\begin{proof}
Fix $p\in(1,\frac{K}{K-1})$. The operator $\cal C$ defined by
$${\cal C}f(z)=\displaystyle\frac{-1}{\pi}\int_{\C}f(w)\left(\frac{1}{w-z}-\frac{1}{w}\right)dA(w)$$
maps functions $f\in L^{2p}(\C)$ to ${\cal C}f\in W^{1,2p}(\C)\subset\Lip_{1-1/p}(\C)$, and
$$\frac{|{\cal C}f(z)-{\cal C}f(w)|}{|z-w|^{1-1/p}}\leq C\,\|f\|_{2p}$$
for any $z,w\in\C$. On the other hand, since $\mu_n,\mu$ are compactly supported on $\D$, there are $L^{2p}(\C)$ functions $h_n, h$ such that
$$\phi_n(z)=z+{\cal C} h_n(z)$$
$$\phi(z)= z +{\cal C}h(z)$$
Actually, $h_n$ and $h$ can be obtained, respectively, as $L^{2p}(\C)$ solutions of the equations $h_n=\B(\mu_nh_n)+\B(\mu_n)$ and $h=\B(\mu\,h)+\B(\mu)$, where $\B$ denotes the Beurling transform, $\B f=p.v\frac{-1}{\pi\,z^2}\ast f$. A computation shows that
$$
(I-\B\mu_n)(\partial(\phi-\phi_n))=\B((\mu-\mu_n)\partial\phi)
$$
Thus, as $I-\B\mu_n$ is a bilipschitz homeomorphism in $L^{2p}(\C)$ whose norm depends only on $\|\mu_n\|\leq\frac{K-1}{K+1}$ and $p$ \cite[Theorem 1]{AIS},
$$\|\partial(\phi-\phi_n)\|_{2p}\leq\|(I-\B\mu_n)^{-1}\|_{2p}\,\|\B((\mu-\mu_n)\,\partial\phi)\|_{2p}$$
and hence $\|\partial(\phi-\phi_n)\|_{2p}\to 0$ as $n\to\infty$. From the Beltrami equation we also get that $\|h-h_n\|_{2p}=\|\overline\partial(\phi-\phi_n)\|_{2p}\to 0$, and hence $J\phi_n$ converges in $L^{2p}$ to $J\phi$. Furthermore, we also have
$$
\aligned
\|\phi-\phi_n\|_{\Lip_{1-1/p}(\C)}&=\|{\cal C}(h-h_n)\|_{\Lip_{1-1/p}(\C)}\leq\|{\cal C}\|\cdot\|h-h_n\|_{2p}
\endaligned
$$
and in particular $\phi_n\to\phi$ uniformly on compact sets.\\
Finally, we note that normalized $K$-quasiconformal mappings are locally H\"older-continuous with exponent $1/K$, with a constant that depends only on $K$. Thus, if $E$ is a compact set, 
$$\aligned
\sup_{z\in E}|\phi^{-1}(z)-\phi_n^{-1}(z)|&=\sup_{z\in E}|\phi_n^{-1}\circ\phi_n\circ\phi^{-1}(z)-\phi_n^{-1}(z)|\\
&\leq C(K)\cdot\sup_{z\in E}|\phi_n\circ\phi^{-1}(z)-z|^{1/K}\\
&= C(K)\cdot\sup_{w\in\phi(E)}|\phi_n(w)-\phi(w)|^{1/K}
\endaligned
$$
Hence $\phi_n^{-1}\to\phi^{-1}$ uniformly on compact sets.
\end{proof}

\begin{theo}\label{zeromeasure}
Let $E\subset\C$ be a compact set, $\alpha\in(0,1]$ and $K\geq 1$. Set $d=\frac{2}{K+1}(1+\alpha K)$. If
$$\H^d(E)=0$$ 
then, $E$ is removable for $K$-quasiregular mappings in $\Lip_\alpha$. 
\end{theo}

\begin{proof}
Assume that $f:\C\rightarrow\C$ is a $\Lip_\alpha(\C)$ mapping, $K$-quasiregular in $\C\setminus E$. We will see that $f$ is in fact $K$-quasiregular in the whole plane $\C$.\\
Let $\mu=\frac{\overline\partial f}{\partial f}$ be the Beltrami coefficient of $f$. Obviously, there is no restriction if we assume that $E\subset\D$ and that $\mu$ is supported on $\D$. Let $\phi$ be the principal solution to $\overline\partial\phi=\mu\,\partial\phi$. Then the function 
$F=f\circ\phi^{-1}$ is holomorphic on $\C\setminus\phi(E)$ and H\"older continuous on $\C$ with exponent $\alpha/K$. In particular, $\overline\partial F$ defines a distribution compactly supported on $\phi(E)$. Indeed, since $F$ is continuous, it is an $L^2_{loc}(\C)$ function, so that its distributional derivatives define continuous linear functionals on compactly supported $W^{1,2}(\C)$ functions. From Weyl's lemma, it is enough to show that $\overline\partial F=0$.\\
For each $\varepsilon>0$ we can consider a finite family of disjoint disks $D_j=D(z_j,r_j)$  such that $\Omega_\varepsilon=\bigcup_j2D_j$ covers $E$ and
$$\sum_{j=1}^n\diam(2D_j)^d\leq \varepsilon$$
We define $\mu_\varepsilon=\mu\,\chi_{\C\setminus\Omega_\varepsilon}$. If $\phi_\varepsilon$ is the principal solution to $\overline\partial\phi_\varepsilon=\mu_\varepsilon\,\partial\phi_\varepsilon$, then by Lemma \ref{convergence} the functions $F_\varepsilon=f\circ\phi_\varepsilon^{-1}$ converge uniformly to $F$ as $\varepsilon\to 0$, since $f\in\Lip_\alpha(\C)$. Thus, for any compactly supported ${\cal C}^\infty$ function $\varphi$,
$$\aligned
\langle\overline\partial F, \varphi\rangle
&=\langle\overline\partial(F-F_\varepsilon),\varphi\rangle+\langle\overline\partial F_\varepsilon,\varphi\rangle\\
&=-\langle(F-F_\varepsilon),\overline\partial\varphi\rangle-\langle\overline\partial F_\varepsilon,\varphi\rangle
\endaligned$$
Here the first term converges to $0$ as $\varepsilon\to 0$. Thus, we just have to care about the second term.\\
Consider a partition of unity $\psi_j$ subordinated to the covering $2D_j$, that is, each $\psi_j$ is a ${\cal C}^\infty$ function, compactly supported on $4D_j$, $|D\psi_j|\leq\frac{C}{\diam(2D_j)}$, and $\sum_{j=1}^n\psi_j=1$ on 
$\Omega_\varepsilon$. We now define $\varphi_j=\psi_j\circ\phi_\varepsilon^{-1}$. Then, each $\varphi_j$ is in $W^{1,2}(\C)$, has compact support in $\phi_\varepsilon(4D_j)$, and $\sum_{j=1}^n\varphi_j=1$ on $\phi_\varepsilon(\Omega_\varepsilon)$. Thus, for any constants $c_j$ we have
\begin{equation}\label{I+II}
\aligned
\langle\overline\partial F_\varepsilon, \varphi\rangle
&=\langle\overline\partial F_\varepsilon, \varphi\sum_{j=1}^n\varphi_j\rangle=\sum_{j=1}^n\langle\overline\partial F_\varepsilon, \varphi\,\varphi_j\rangle=\sum_{j=1}^n\langle\overline\partial (F_\varepsilon-c_j), \varphi\,\varphi_j\rangle\\
&=-\sum_{j=1}^n\langle (F_\varepsilon-c_j), \overline\partial\varphi\,\varphi_j\rangle-\sum_{j=1}^n\langle (F_\varepsilon-c_j), \varphi\,\overline\partial\varphi_j\rangle
=I+II
\endaligned
\end{equation}
and now we have just to bound both sums independently.\\
On one hand, 
$$
\aligned
\left|\sum_{j=1}^n\langle (F_\varepsilon-c_j), \overline\partial\varphi\,\varphi_j\rangle\right|
&=\left|\sum_{j=1}^n\int_{\phi_\varepsilon(4D_j)}(F_\varepsilon(z)-c_j)\,\overline\partial\varphi(z)\,\varphi_j(z)\,dA(z)\right|\\
&\leq\sum_{j=1}^n\int_{\phi_\varepsilon(4D_j)}|F_\varepsilon(z)-c_j|\,|\overline\partial\varphi(z)|\,|\varphi_j(z)|\,dA(z)\\
&\leq\|D\varphi\|_\infty\,\sum_{j=1}^n\int_{4D_j}|f(w)-c_j|\,J\phi_\varepsilon(w)dA(w)
\endaligned
$$
Now observe \cite{A} that for any $p\in(1,\frac{K}{K-1})$, $J\phi_\varepsilon\in L^p_{loc}(\C)$. Using this, together with the fact that $f\in\Lip_\alpha$, we obtain for $I$ the bound
$$
\aligned
|I|
&\leq\,C\,\|D\varphi\|_\infty\,\sum_{j=1}^n\left(\int_{4D_j}J\phi_\varepsilon(w)^pdA(w)\right)^{1/p}\diam(4D_j)^{\alpha+2/q}\\
&\leq\,C\,\|D\varphi\|_\infty\,\left(\sum_{j=1}^n\int_{4D_j}J\phi_\varepsilon(w)^pdA(w)\right)^{1/p}\,\left(\sum_{j=1}^n\diam(4D_j)^{q\alpha+2}\right)^{1/q}\\
\endaligned
$$
where $q=\frac{p}{p-1}$. But both sums converge to $0$ as $\varepsilon$ tends to $0$, so that we get $I=0$. On the other hand, an analogous computation for $II$ shows that 
$$
\aligned
\left|\sum_{j=1}^n\langle (F_\varepsilon-c_j),\varphi\,\overline\partial\varphi_j\rangle\right|
&=\left|\sum_{j=1}^n\int_{\phi_\varepsilon(4D_j)}(F_\varepsilon(z)-c_j)\,\varphi(z)\,\overline\partial\varphi_j(z)\,dA(z)\right|\\
&\leq\sum_{j=1}^n\int_{\phi_\varepsilon(4D_j)}|F_\varepsilon(z)-c_j|\,|\varphi(z)|\,|\overline\partial\varphi_j(z)|\,dA(z)\\
&\leq\|\varphi\|_\infty\,\sum_{j=1}^n\int_{4D_j}|f(w)-c_j|\,|\overline\partial\varphi_j(\phi_\varepsilon(w))|\,J\phi_\varepsilon(w)dA(w)
\endaligned
$$
The chain rule gives
$$|D\varphi_j(\phi_\varepsilon)|\leq|D\psi_j|\,|D\phi_\varepsilon^{-1}(\phi_\varepsilon)|$$
Thus, we have
$$\aligned
|II|
&\leq\|\varphi\|_\infty\,\sum_{j=1}^n\left(\int_{4D_j}|D\phi_\varepsilon^{-1}(\phi_\varepsilon(w))|\,J\phi_\varepsilon(w)\,dA(w)\right)\,\diam(4D_j)^{\alpha-1}\\
&=\|\varphi\|_\infty\,\sum_{j=1}^n\left(\int_{\phi_\varepsilon(4D_j)}|D\phi_\varepsilon^{-1}(z)|\,dA(z)\right)\,\diam(4D_j)^{\alpha-1}\\
&\leq\|\varphi\|_\infty\,\sum_{j=1}^n\left(\int_{4D_j}|D\phi_\varepsilon(w)|\,dA(w)\right)\,\diam(4D_j)^{\alpha-1}
\endaligned
$$
But $J\phi_\varepsilon^{1/2}$ defines a doubling measure, so that
$$\aligned
\int_{4D_j}|D\phi_\varepsilon(w)|\,dA(w)\leq C\,\int_{2D_j}J\phi_\varepsilon(z)^{1/2}dA(z)
\endaligned$$
where $C$ depends only on $K$. On the other hand, $\phi_\varepsilon$ is conformal in $\Omega=\cup_j2D_j$. Hence, from the improved integrability results for $K$-quasiconformal mappings \cite{AN}, we have $J\phi_\varepsilon\in 
L^\frac{K}{K-1}(\Omega_\varepsilon)$ and
$$\int_{\Omega_\varepsilon}J\phi_\varepsilon(w)^\frac{K}{K-1}\,dA(w)\leq 1$$
Hence, if we choose $p=\frac{2K}{K-1}$ and $q=\frac{2K}{K+1}$, then $q(\alpha-1)+2=d$ and we get
$$
\aligned
|II|
&\leq\,C\,\|\varphi\|_\infty\,\sum_{j=1}^n\left(\int_{2D_j}J\phi_\varepsilon(w)^{1/2}dA(w)\right)\,\diam(2D_j)^{\alpha-1}\\
&\leq\,C\,\|\varphi\|_\infty\,\sum_{j=1}^n\left(\int_{2D_j}J\phi_\varepsilon(w)^{p/2}dA(w)\right)^{1/p}\,\diam(2D_j)^{\alpha-1+2/q}\\
&\leq\,C\,\|\varphi\|_\infty\,K^{1/2}\,\left(\int_{\Omega_\varepsilon}J\phi_\varepsilon(w)^{p/2}dA(w)\right)^{1/p}\,\left(\sum_{j=1}^n\diam(2D_j)^{q(\alpha-1)+2}\right)^{1/q}
\endaligned
$$
But here the first sum remains bounded, while the second converges to $0$ as $\varepsilon\to 0$.
\end{proof}

We wish to emphasize that the improved integrability properties for $K$-quasiconformal mappings \cite{AN} is precisely what allows us to jump from sets of dimension strictly smaller than $d$ to sets of null $d$-dimensional Hausdorff measure. Notice that in the particular case $\alpha=1$, if $\H^2(E)>0$ then there exists a nonconstant $\Lip_1(\C)$ function $f$ which is holomorphic (and thus $K$-quasiregular, for any $K$) on $\C\setminus E$ (e.g.\cite{U}). As we said in the introduction, it is also well known that compact sets $E$ with $\H^{1+\alpha}(E)>0$ are not removable for analytic functions in $\Lip_\alpha$, for any $\alpha\in (0,1)$.\\
\\
Note that Theorem \ref{zeromeasure} may be stated in a more general sense. Assume that $f$ is $K$-quasiregular in $\C\setminus E$ and, instead of equation (\ref{lipalfa}), suppose that
$$|f(z)-f(w)|\leq \omega(|z-w|)$$
for any $z,w\in\C$ with  $|z-w|<1$, where $\omega(r)=\omega_f(r)$ is any bounded continuous non-decreasing function $\omega:[0,1)\to[0,\infty)$ such that $\omega(0)=0$. If $\H^h(E)=0$ for $h(t)=t^\frac{2}{K+1}\,\omega(2t)^\frac{2K}{K+1}$, then $f$ is in fact $K$-quasiregular on the whole plane. As particular cases, $\omega(r)=C\,r^\alpha$ gives the previous Theorem, while if $f\in\lip_\alpha(\C)$, that is, if $\frac{\omega_f(r)}{r^\alpha}\to 0$ as $r\to 0$, then $\M^d_\ast(E)=0$ is enough (here $\M^d_\ast$ denotes the lower $d$-dimensional Hausdorff content). Finally, if we only assume $f$ to be continuous, then our argument shows that the condition $\M^\frac{2}{K+1}_\ast(E)=0$ suffices.\\
\\
The same procedure also gives removability results in the more general setting of finite distortion mappings. More precisely, subexponentially integrable distortion mappings are locally in $W^{1,p}$ for any $p<2$. In this context, we obtain the following.

\begin{coro}
Suppose that $f\in\Lip_\alpha(\C)$, and let $E\subset\D$ be a compact set. Assume that $f$ is a mapping of finite distortion on $\C\setminus E$, with distortion function $K$ such that $K(z)=1$ whenever $|z|>1$, and
\begin{equation}\label{subexponential}
\exp\left(\frac{K}{\log(e+K)}\right)\in L^p(\D)
\end{equation}
for some $p>0$. If $\dim(E)<2\alpha$, then $f$ is a mapping of finite distortion on $\C$.
\end{coro}

Actually, here one can replace condition (\ref{subexponential}) by the more general condition given by
$$\exp\left(\A(K)\right)\in L^p(\D)
$$
where $\A:[1,\infty)\rightarrow[1,\infty)$ is any smoothly increasing function, such that $\lim_{t\to\infty}t\,\A'(t)=\infty$, $t\mapsto e^{\A(t)}$ is convex, and $\int_1^\infty\frac{\A(t)}{t^2}dt=\infty$ \cite{IKO}. On the other hand, if we ask $K$ to be in a stronger space, then one can use the sharp integrability results in \cite{FKZ} and obtain bigger removable singularities.

\begin{coro}
There exists a positive constant $p_0>0$ with the following property. Let $f\in\Lip_\alpha(\C)$, and $E\subset\D$ be a compact set. Assume that $f$ is a mapping of finite distortion on $\C\setminus E$, with distortion function $K$ such that $K(z)=1$ whenever $|z|>1$, and
$$\exp\left(K\right)\in L^p(\D)$$
for some $p>p_0$. If $\H^{2\alpha}(E)$ is $\sigma$-finite, then $f$ is a mapping of finite distortion on $\C$.
\end{coro}

{\bf{Acknowledgements}}. I am specially grateful to my advisors J. Mateu and J. Orobitg for interessing discussions. I am also grateful to K. Astala for interesting conversations, and for his suggestion about the critical index.

\vskip 1cm
\begin{itemize}

\item[]{Departament de Matem\`atiques, Facultat de Ci\`encies\\Campus de la Universitat Aut\`onoma de Barcelona\\ 08193-Bellaterra, Barcelona (Spain)\\
 {albertcp@mat.uab.es}\\ Tel. +34 93 581 45 45\\Fax +34 93 581 27 90}

\end{itemize}

\end{document}